\documentclass[reqno,11pt,reqnopdfmx]{amsart}

\usepackage[top=25truemm,bottom=25truemm,left=25truemm,right=25truemm]{geometry}
\usepackage[bookmarksnumbered,colorlinks,plainpages]{hyperref}
\usepackage{amsmath,amsthm,amssymb,amsfonts,amscd,mathtools,color}
\usepackage[all]{xy}
\usepackage{graphicx} 
\numberwithin{equation}{section}

\theoremstyle{plain}
	\newtheorem{thm}{Theorem}[section]
	\newtheorem*{thm*}{Theorem}
	\newtheorem{cor}[thm]{Corollary}
	\newtheorem{lem}[thm]{Lemma}
	
	\newtheorem{prop}[thm]{Proposition}
	
	\newtheorem*{conj*}{Conjecture}
	
\theoremstyle{definition}
	\newtheorem{defn}[thm]{Definition}

\theoremstyle{remark}
	\newtheorem{rem}[thm]{Remark}
	
	\newtheorem*{pf}{Proof}
\def\CC{{\mathbb C}}

\def\HH{{\mathbb H}}

\def\PP{{\mathbb P}}
\def\QQ{{\mathbb Q}}
\def\RR{{\mathbb R}}

\def\ZZ{{\mathbb Z}}

\def\D{{\mathcal D}}
\def\E{{\mathcal E}}
\def\F{{\mathcal F}}

\def\L{{\mathcal L}}

\def\O{{\mathcal O}}

\def\R{{\mathcal R}}
\def\S{{\mathcal S}}

\def\U{{\mathcal U}}

\def\Aut{{\rm Aut}}

\def\coh{{\rm coh}}

\def\Exc{{\rm Exc}}
\def\Ext{{\rm Ext}}
\def\FEC{{\rm FEC}}

\def\Hom{{\rm Hom}}
\def\Ind{{\rm Ind}}

\def\mod{{\rm mod}}

\def\vd{\vec{\Delta}}
\begin{document}
\title{Number of full exceptional collections modulo spherical twists for elliptic orbifolds}

\author{Atsushi Takahashi}
\address{Department of Mathematics, Graduate School of Science, Osaka University,
Toyonaka Osaka, 560-0043, Japan}
\email{takahashi@math.sci.osaka-u.ac.jp}
\author{Hongxia Zhang}
\address{School of Mathematical Sciences, Xiamen University,
Xiamen, 361005, Fujian, China.}
\email{hongxiazhang@stu.xmu.edu.cn}

\maketitle
\thispagestyle{empty}
\begin{abstract}
This paper calculates the number of full exceptional collections modulo an action of a group as the set generated by spherical twists for an abelian category of coherent sheaves on an orbifold projective line with a zero orbifold Euler characteristic. This is done by a recursive formula naturally generalizing the one for the Dynkin case by Deligne whose categorical interpretation is due to Obaid-Nauman-Shammakh-Fakieh-Ringel and an abelian category of coherent sheaves on an orbifold projective line with a positive orbifold Euler characteristic is due to Otani-Shiraishi-Takahashi.
\end{abstract}
\section{Introduction}

Inspired by the correspondence between maximal chains in the proof of noncrossing partitions for a Dynkin quiver $\vec{\triangle}$ and complete exceptional sequences in the derived category $\D^b(\CC\vec{\Delta})$ of finitely generated modules over the path algebra $\CC\vec{\Delta}$, Obaid-Nauman-Shammakh-Fakien-Ringel {\cite{ONSFR}} found the following recursive formula for $e(\D^b(\CC\vec{\Delta})) :=|\FEC(\D^b(\CC\vec{\Delta})) / \ZZ^\mu|$, where $\FEC(\D^b(\CC\vec{\Delta}))$ is the set of isomorphism classes of full exceptional collection in $\D^b(\CC\vec{\Delta})$, $\mu$ is the number of vertices of $\vec{\Delta}$ and $\ZZ^\mu$ is the group whose $i$-th generator acts by the translation functor $[1]$ on the $i$-th object of a full exceptional collection. $$e(\D^b(\CC\vec{\Delta}))=\frac{h}{2}\sum_{\nu\in\Delta_0}e(\D^b(\CC\vec{\Delta}^{(\nu)})),$$
where $\vec{\Delta}^{(\nu)}$ is the full subquiver of $\vec{\Delta}$ given by removing the vertex $\nu$ and arrows connecting with $\nu$ and $h\in \ZZ_{\geq 0}$ is the Coxeter number. As is mentioned in {\cite{ONSFR}}, the above recursive formula is a categorification of Deligne's recursive formula {\cite{De}} which is originally answered to Looijenga's conjecture for simple singularities concerning the number of distinguished bases of vanishing cycles modulo signs and the degree of a map $LL$, called the Lyashko-Looijenga map, describing the topology of the bifurcation set {\cite{Lo}}. As a summary of the above story, the following equation holds:
$$e(\D^b(\CC\vec{\Delta}))= \frac{\mu!}{d_1\cdots d_\mu}h^\mu={\rm deg}LL,$$
where $2\le d_1\le\dots\le d_\mu=h$ are degrees of algebraic independent invariants of the Weyl group for the Dynkin diagram $\Delta$. Here $LL$ was originally defined by using the characteristic polynomial associated to the discriminant for the singularity, or equivalently, the square of Jacobian for the Wely group invariant theory. More generally, $LL$ can be defined for massive $F$-manifolds (see {\cite{DZ, He}}), In {\cite{HR}}, Hertling-Roucairol extended the above correspondence between the number of distinguished bases of vanishing cycles modulo signs and the degree of the Lyashko-Looijenga map to simple elliptic singularities with Legendre normal forms. Therefore, it is very important to understand what happens when the setting falls between the two classes, corresponding to simple singularities and the simple elliptic singularities.

Recently, {\cite{OST}} generalized the above recursive formula and equality for the derived category $\D^b(\PP_{A}^1)$ of coherent sheaves over an orbifold projective line $\PP_{A}^1$ with three orbifold points $(0,1, \infty)$ whose orders are given by $A:=(a_1, a_2, a_3)$ satisfying $\chi_A:=\frac{1}{a_1}+\frac{1}{a_2}+\frac{1}{a_3}-1\textgreater 0,$ that is, orbifold projective line of domestic type. They define $e(\D^b(\PP_{A}^1))$ from the perspective of mirror symmetry and the space of stability conditions. They reach to definition of $e(\D^b(\PP_{A}^1))$ as
$$e(\D^b(\PP_{A}^1)):=|\FEC(\D^b(\PP_{A}^1))/\langle {\rm ST}(\D^b(\PP_{A}^1)) , \ZZ^{\mu_A}\rangle |.$$
where ${\rm ST}(\D^b(\PP_{A}^1))$ generated by the spherical twist $=(-)\otimes\O(\vec{c})$, $\mu_A=a_1+a_2+a_3-1.$
They prove the following holds
$$e(\D^b(\PP_{A}^1)) =\frac{\mu_A!}{a_1!a_2! a_3!\chi_{A}}a_1^{a_1}a_2^{a_2} a_3^{a_3}.$$

The purpose of this paper is to consider $e(\D^b(\PP_{A,\Lambda}^1))$ for the derived category $\D^b(\PP_{A,\Lambda}^1)$ of coherent sheaves over an orbifold projective line $\PP_{A,\Lambda}^1$ with orbifold points $\Lambda:=(\lambda_1,\dots, \lambda_r)$, whose orders are given by $(a_1,\dots, a_r)$ satisfying $\chi_A:=\frac{1}{a_1}+\cdots+\frac{1}{a_r}-1=0,$ where
$\Lambda:=(\lambda_1,\dots, \lambda_r)$ the ordered set of pairwise distinct points on $\PP^1$ normalized as
$\lambda_1:=\infty$, $\lambda_2:=0$ and $\lambda_3:=1$.

Let $\Aut(\D^b(\PP^1_{A,\Lambda}))$ be the group of autoequivalences of derived coherent category $\D^b(\PP^1_{A,\Lambda}),$ consider the quotient of the subgroup ${\rm ST}(\D^b(\PP^1_{A,\Lambda}))$ of $\Aut(\D^b(\PP^1_{A,\Lambda}))$ by ${\rm ST}(\D^b(\PP^1_{A,\Lambda}))\cap \ZZ[1]$
and denote it by $\overline{{\rm ST}(\D^b(\PP^1_{A,\Lambda}))}$.

Let
$\Gamma:={\rm SL}(2;\ZZ)$ be the special linear group of $2\times 2$ matrices with integer coefficients and determinant $1$.
Denote by $\overline{\Gamma}$ the quotient of $\Gamma$ by its center. Let
$$\Gamma(l)=\left\{\left.\begin{pmatrix} a & b\\ c & d\end{pmatrix}\, \right|\, a,d\equiv 1\ ({\rm mod}\ l),\ b,c\equiv 0 \ ({\rm mod}\ l)\right\}.$$
whose image in $\overline{\Gamma}$ is denoted by $\overline{\Gamma(l)}$. Denote by $\ell_A:={\rm lcm}(a_1,\dots, a_r)$. The first main result is the group isomorphism between ${\rm ST}(\D^b(\PP^1_{A,\Lambda}))/\ZZ[2]$ and $\Gamma(\ell_A)$.

\begin{thm}[Theorem~\ref{thm:main1}]
Suppose $\chi_A=0$. The group ${\rm ST}(\D^b(\PP^1_{A,\Lambda}))$ is a subgroup of ${\rm Br}_3$ satisfying
$${\rm ST}(\D^b(\PP^1_{A,\Lambda}))/\ZZ[2]\cong \Gamma(\ell_A),\quad \overline{{\rm ST}(\D^b(\PP^1_{A,\Lambda}))}\cong \overline{\Gamma(\ell_A)}.$$
\end{thm}
The second main result is the equality $e(\D^b(\PP_{A,\Lambda}^1))$. We obtain the following recursion relation:
$$e(\D^b(\PP_{A,\Lambda}^1))=\sum_{i=1}^r\sum_{j=1}^{a_i-1}\binom{\mu_A-1}{j-1}\cdot e(\D^b(\PP^1_{A_{i,j},\Lambda}))\cdot e(\D^b(\CC\vec{A}_{j-1}))\cdot a_i\cdot \frac{[\overline{\Gamma}:\overline{\Gamma(\ell_A)}]}{\ell_A},$$
where $A_{i,j}:=(a_1,\dots a_{i-1},a_i-j,a_{i+1}\dots ,a_r)$. It turns out that $\PP^1_{A_{i,j},\Lambda}$ is an orbifold projective line of domestic type, according to \cite[Theorem 4.10]{OST}, it follows that $$e(\D^b(\PP^1_{A_{i,j},\Lambda}))=\frac{(\mu_A-j)!}{(a_1-j)!a_2!\cdots a_r!\chi_{A_{1,j}}}(a_1-j)^{a_1-j}a_2^{a_2}\cdots a_r^{a_r}.$$
Furthermore, we obtain the following

\begin{thm}[Theorem~\ref{thm:main2}]
Suppose $\chi_A=0$.
Then we have
$$e(\D^b(\PP_{A,\Lambda}^1))=\frac{\mu_A!}{a_1!\cdots a_r!}a_1^{a_1}\cdots a_r^{a_r}\cdot
\frac{\displaystyle\sum_{i=1}^r a_i^2 (a_i-1)}{2\mu_A}\cdot \frac{[\overline{\Gamma}:\overline{\Gamma(\ell_A)}]}{\ell_A}.$$
In particular,
\[
|{\rm FEC}(\D^b(\PP_{A,\Lambda}^1))/\overline{\Gamma}|=
\frac{e(\D^b(\PP_{A,\Lambda}^1))}{[\overline{\Gamma}:\overline{\Gamma(\ell_A)}]}
=
\frac{\mu_A!}{a_1!\cdots a_r!}a_1^{a_1}\cdots a_r^{a_r}\cdot
\frac{\displaystyle\sum_{i=1}^r a_i^2 (a_i-1)}{2\ell_A\mu_A}.
\]
\end{thm}
As a consequence (Corollary~\ref{cor:3.8}), we have another expression of the degree ${\rm deg}LL^{alg}$ of the Lyashko-Looijienga map
for the universal unfolding of a simple elliptic singularity in the Legendre normal form calculated by Hertling-Roucairol~\cite{HR}.
More concretely, we have
\[
{\rm deg}LL^{alg}=|{\rm FEC}(\D^b(\PP_{A,\Lambda}^1))/\overline{\Gamma(2)}|=\frac{\mu_A!}{a_1!\cdots a_r!}a_1^{a_1}\cdots a_r^{a_r}\cdot
\frac{\displaystyle\sum_{i=1}^r a_i^2 (a_i-1)}{2\ell_A\mu_A}\cdot 6,
\]
In particular, it is important to note that ours is given in terms of $A$ by the above  recursion relation which can be considered as ``the cut-and-join equation",
while their formula
\[
{\rm deg}LL^{alg}= \frac{\mu!\cdot \frac{1}{2}\sum_{j=2}^{\mu-1}\frac{1}{{\rm deg_\textbf{w}}t_j}} {\prod_{j=2}^{\mu-1}{{\rm deg}_\textbf{w}t_j}},
\]
can be considered as a variation of ELSV formulas, which is based on weights ${\rm deg}_\textbf{w}t_i$ of homogeneous coordinates $t_i$ on the base space of the universal unfolding given by Table~\ref{tab: weights}.
\begin{table}[h]
\begin{tabular}{|c||c|c|c|}
\hline
Singularity Type  & Universal Unfolding  & ${\rm deg}_\textbf{w}t_1,\dots, {\rm deg}_\textbf{w}t_\mu$\\ \hline \hline
$\widetilde{E}_6$ & $z_2(z_2-z_1)(z_2-\lambda z_1)-z_1z_3^2+\sum_{i=1}^6 t_i\phi_i(z)$ & $1,\frac{2}{3},\frac{2}{3},\frac{2}{3},\frac{1}{3},\frac{1}{3},\frac{1}{3},0$ \\ \hline
$\widetilde{E}_7$ & $z_1z_2(z_2-z_1)(z_2-\lambda z_1)+z_3^2+\sum_{i=1}^7 t_i\phi_i(z)$ & $1,\frac{3}{4},\frac{3}{4},\frac{2}{4},\frac{2}{4},\frac{2}{4},\frac{1}{4},\frac{1}{4},0$\\ \hline
$\widetilde{E}_8$ & $z_2(z_2-z_1^2)(z_2-\lambda z_1^2)+z_3^2+\sum_{i=1}^8 t_i\phi_i(z)$ & $1,\frac{5}{6},\frac{4}{6},\frac{4}{6},\frac{3}{6},\frac{3}{6},\frac{2}{6},\frac{2}{6},\frac{1}{6},0$\\ \hline
\end{tabular}
\vspace{5pt}
\caption{Weights of homogeneous coordinates.}\label{tab: weights}
\end{table}

{\bf Acknowledgements.}
The first named author is grateful to Claus Hertling for valuable discussions.
This work has been supported by JSPS KAKENHI Grant Number JP21H04994. The second named author is also supported by China Scholarship Council. The authors are grateful to Yuuki Shiraishi for valuable comments.

\vspace{5mm}

{\bf Notation.}
Throughout this paper, the translation functor of a triangulated category will be denoted by $[1]$. For a finite dimensional $\CC$-algebra $A$, the bounded derived category of finitely generated right $A$-modules is denoted by $\D^b(A):=\D^b\mod(A)$. We also denote by $\D^b(X)$ the derived category $\D^b\coh(X)$ of coherent sheaves on an orbifold $X$. For a triangulated category $\D$, the group of autoequivalences of $\D$ is denoted by $\Aut(\D)$. 
\section{Preliminaries}\label{sec : preliminaries}

\subsection{Set of full exceptional collections} First, we recall the notion of a full exceptional collection and basic properties.
\begin{defn}
Let $\D$ be a $\CC$-linear triangulated category.
\begin{enumerate}
\item
An object $E \in \D$ is called {\em exceptional} if $\Hom_\D(E,E) \cong \CC$ and $\Hom_\D(E,E[p]) \cong 0$ when $p \ne 0$.
\item
For an exceptional object $E\in\D$, $E^\perp$ is a full triangulated subcategory of $\D$ defined as
\begin{equation*}
E^\perp := \{ X \in \D \mid \Hom_\D(E, X[p]) = 0 ~ \text{for all} ~ p \in \ZZ \},
\end{equation*}
\item
An ordered set $\E = (E_1, \dots, E_\mu)$ consisting of exceptional objects $E_1, \dots, E_\mu$ is called {\em exceptional collection} if $\Hom_\D(E_i, E_j[p]) \cong 0$ for all $p \in \ZZ$ and $i > j$.
\item
An exceptional collection $\E$ is called {\em full} if the smallest full triangulated subcategory of $\D$ containing all elements in $\E$ is equivalent to $\D$ as a triangulated category.
\item
Two full exceptional collections $\E=(E_1, \dots, E_\mu)$, $\F=(F_1, \dots, F_\mu)$ in $\D$ are said to be {\em isomorphic} if $E_i \cong F_i$ for all $i = 1, \dots, \mu$.
The set of isomorphism classes of full exceptional collections in $\D$ is denoted by $\FEC(\D)$.
\end{enumerate}
\end{defn}

An object $E \in \D$ is said to be {\em indecomposable} if $E$ is nonzero and $E$ has no direct sum decomposition $E\cong A\oplus B,$ where $A$ and $B$ are nonzero objects in $\D$. Note that an exceptional object is indecomposable.

\begin{defn}
Let $(E, F)$ be an exceptional collection. Define two objects $\R_F E$ and $\L_E F$ by the following exact triangles respectively:
\begin{equation*}
\R_F E \longrightarrow E \overset{{\rm ev}^*}{\longrightarrow} \Hom^\bullet_\D(E, F)^* \otimes F\longrightarrow(\R_F E)[1],
\end{equation*}
\begin{equation*}
(\L_E F)[-1]\longrightarrow\Hom^\bullet_\D(E, F) \otimes E \overset{\rm ev}{\longrightarrow} F \longrightarrow \L_E F,
\end{equation*}
where $(-)^*$ denotes the duality $\Hom_\CC(-, \CC)$. The object $\R_F E$ (\emph{resp.} $\L_E F$) is called the \emph{right mutation} of $E$ through $F$ (\emph{resp.} left mutation of $F$ through $E$). Then, $(F, \R_F E)$ and $(\L_E F, E)$ form new exceptional collections.
\end{defn}

\begin{rem}
Our definition of mutations differs from the usual one (cf. \cite[Section 1]{BP}. In our notation, the usual left and right mutations are given by $\R_F E[1]$ and $\L_E F[-1]$, respectively.
\end{rem}

The Artin's {\em braid group} ${\rm Br}_{\mu}$ on $\mu$-stands is a group presented by the following generators and relations:
\begin{description}
\item[{\bf Generators}] $\{b_i~|~i=1,\dots, \mu-1\}$
\item[{\bf Relations}] $b_{i}b_{j}=b_{j}b_{i}$ for $|i-j|\ge 2$, $b_{i}b_{i+1}b_{i}=b_{i+1}b_{i}b_{i+1}$ for $i=1,\dots, \mu-2$.
\end{description}
Consider the group ${\rm Br}_{\mu} \ltimes \ZZ^{\mu}$, the semi-direct product of the braid group ${\rm Br}_{\mu}$ and the
abelian group $\ZZ^{\mu}$, defined by the group homomorphism ${\rm Br}_{\mu} \longrightarrow {\mathfrak S}_{\mu} \longrightarrow \Aut_{\ZZ}(\ZZ^{\mu})$, where the first homomorphism is $b_{i}\mapsto (i, i+1)$ and the second one is induced by the natural actions of the symmetric group ${\mathfrak S}_{\mu}$ on $\ZZ^{\mu}$.
\begin{prop}[{cf. Bondal--Polishchuk~\cite[Proposition 2.1]{BP}}]\label{prop : braid}
Let $\D$ be a $\CC$-linear triangulated category.
The group ${\rm Br}_{\mu} \ltimes \ZZ^{\mu}$ acts on $\FEC(\D)$ by mutations and transformations:
\begin{gather*}
b_{i} \cdot (E_{1}, \dots, E_{\mu}) \coloneqq (E_{1}, \dots, E_{i-1}, E_{i+1}, \R_{E_{i+1}} E_i, E_{i+2}, \dots, E_{\mu}), \\
b^{-1}_{i} \cdot (E_{1}, \dots, E_{\mu}) \coloneqq (E_{1}, \dots, E_{i-1}, \L_{E_{i}} E_{i+1}, E_{i}, E_{i+2},\dots, E_{\mu}), \\
e_{i} \cdot (E_{1}, \dots, E_{\mu}) \coloneqq (E_{1}, \dots, E_{i-1}, E_{i} [1], E_{i+1}, \dots, E_{\mu}),
\end{gather*}
where $e_{i}$ is the $i$-th generator of $\ZZ^{\mu}$.
\qed
\end{prop}

In this paper, we shall consider the derived category $\D^b(\mathcal{A})$ of a hereditary abelian category $\mathcal{A}$, namely, an abelian category $\mathcal{A}$ satisfying $\Ext^p_{\mathcal{A}}(X,Y)=0$ for any $p\geq 2$ and objects $X, Y\in \mathcal{A}.$ Indecomposable objects and exceptional objects in such derived category has the property that

$$\Ind(\D^b(\mathcal{A}))\cong \bigsqcup_{p\in \ZZ}\Ind(\mathcal{A})[p],\;\;\; \Exc(\D^b(\mathcal{A}))\cong \bigsqcup_{p\in \ZZ}\Exc(\mathcal{A})[p],$$
where $\Ind(\mathcal{A})[p]$ (\emph{resp.} $\Exc(\mathcal{A})[p]$) denotes the set of isomorphism classes of indecomposable (\emph{resp.} exceptional) objects $E\in \mathcal{A}$ translated by $p$.

Note that the group of autoequivalences $\Aut(\D)$ of $\D$ also acts on $\FEC(\D)$ by
\[
\Phi(E_1, \dots, E_\mu) \coloneqq (\Phi(E_1), \dots, \Phi(E_\mu)), \quad \Phi \in \Aut(\D).
\]

\subsection{Orbifold projective lines}
For $A=(a_1,\dots, a_r)$, $a_i\in\ZZ_{\ge 2}$, set
\[
\ell_A:={\rm lcm}(a_1,\dots, a_r),\quad
\mu_A:=2+\sum_{i=1}^r\left(a_i-1\right),\quad
\chi_A:=2+\sum_{i=1}^r\left(\frac{1}{a_i}-1\right).
\]
Note that $\chi_A=0$ if and only if $A=(2,2,2,2)$, $(3,3,3)$, $(2,4,4)$ or $(2,3,6)$.
For the later convenience, we give a table of $\ell_A$ and $\mu_A$ for these cases:
\begin{table}[h]
\begin{tabular}{|c||c|c|c|c|}
\hline
$A$ & $(2,2,2,2)$ & $(3,3,3)$ & $(2,4,4)$ & $(2,3,6)$ \\ \hline \hline
$\ell_A$ & $2$ & $3$ & $4$ & $6$ \\ \hline
$\mu_A$ & $6$ & $8$ & $9$ & $10$\\ \hline
\end{tabular}
\vspace{5pt}
\caption{Classification of $A$ with $\chi_A=0$.}\label{tab: classification A}
\end{table}
\begin{defn}[{Geigle--Lenzing~\cite[Section 1.1]{GL}}]
Let $A=(a_1,\dots, a_r)$, $a_i\in\ZZ_{\ge 2}$ and $\Lambda:=(\lambda_1,\dots, \lambda_r) $ the ordered set of pairwise distinct points on $\PP^1$ normalized as
$\lambda_1:=\infty$, $\lambda_2:=0$ and $\lambda_3:=1$.
\begin{enumerate}
\item Denote by $L_A$ an abelian group generated by $r$-letters $\vec{x_i}$, $i = 1, \dots, r$ defined as the quotient
\begin{equation*}
L_A \coloneqq \left. \bigoplus_{i = 1}^r \ZZ \vec{x}_i \middle/ \big\langle a_i \vec{x}_i - a_j \vec{x}_j \mid i, j = 1, \dots, r \big\rangle \right. .
\end{equation*}
In $L_A$ we have $a_1 \vec{x}_1 = \cdots = a_r \vec{x}_r$, which will be denoted by $\vec{c}_A$.
\item Define a $\CC$-algebra $S_{A,\Lambda}$ as the quotient of $\CC[x_1, \dots, x_r]$ by the ideal generated by
\begin{equation*}
x_i^{a_i} - x_2^{a_2} + \lambda_i x_1^{a_1},\quad i=3,\dots, r.
\end{equation*}
\end{enumerate}
\end{defn}
Note that $L_A$ is an abelian group of rank one and $S_{A,\Lambda}$ is naturally $L_A$-graded.

\begin{defn}[{cf.~\cite[Section 1.1]{GL}}]
The {\em orbifold projective line} $\PP^1_{A,\Lambda}$ of type $(A,\Lambda)$ is a quotient stack defined by
\begin{equation*}
\PP_{A,\Lambda}^1 \coloneqq \left[ \left( {\rm Spec} (S_{A,\Lambda}) \backslash \{ 0 \} \right) / {\rm Spec} ({\CC L_A}) \right].
\end{equation*}
\end{defn}
The orbifold projective line is a Deligne--Mumford stack whose coarse moduli space is a smooth projective line $\PP^1$,
which can be considered as a smooth $\PP^1$ with orbifold points $\lambda_1,\dots, \lambda_r$ of orders $a_1,\dots, a_r$ .

Denote by ${\rm gr}^{L_A}(S_{A,\Lambda})$ (resp. ${\rm gr}^{L_A}_0(S_{A,\Lambda})$) the abelian category of
finitely generated (resp. finite length) $L_A$-graded $S_{A,\Lambda}$-modules.
The abelian category $\coh(\PP^{1}_{A,\Lambda})$ of coherent sheaves on $\PP^{1}_{A,\Lambda}$ is given by
\[
\coh(\PP^{1}_{A,\Lambda}) = {\rm gr}^{L_A}(S_{A,\Lambda}) / {\rm gr}^{L_A}_0(S_{A,\Lambda}),
\]
which is hereditary (see \cite[Section~1.8 and Section~2.2]{GL}).

For $\vec{l} \in L_{A}$, define a sheaf $\O (\vec{l})$  by
\begin{equation*}
\O (\vec{l}) \coloneqq [ S_{A,\Lambda} (\vec{l}) ] \in \coh(\PP_{A,\Lambda}^1)
\end{equation*}
where $(S_{A,\Lambda} (\vec{l}))_{\vec{l}^{\prime}}:=(S_{A,\Lambda} )_{\vec{l}+\vec{l}^{\prime}}$.

\begin{defn}[{\cite[Section 2.5]{GL}}]
For an ordinary point $\lambda \in \PP^{1} \setminus \{ \lambda_1, \dots, \lambda_r \}$, define $S_\lambda$ by
\begin{equation}\label{eq:ordinary simple}
0\rightarrow \O \xrightarrow{x_2^{a_{2}} - \lambda x_1^{a_{1}}} \O (\vec{c}) \rightarrow S_\lambda \rightarrow 0.
\end{equation}
For an orbifold point $\lambda_i$ and $k=0,\dots, a_i-1$, define $S_{i;k}$ by
\begin{equation}
0\rightarrow \O ((k-1) \vec{x}_{i}) \xrightarrow{x_i} \O (k\vec{x}_{i}) \rightarrow S_{i;k} \rightarrow 0.
\end{equation}
\end{defn}

Each coherent sheaf on $\PP_{A,\Lambda}^1$ is decomposed as a direct sum of a vector bundle and a torsion sheaf.
Denote by ${\rm vect}(\PP_{A,\Lambda}^1)$ (resp. $\coh_0 (\PP_{A,\Lambda}^1)$) the full subcategory of $\coh(\PP_{A,\Lambda}^1)$
consisting of vector bundles (resp. torsion sheaves).
It is known that the abelian category $\coh_0 (\PP_{A,\Lambda}^1)$ decomposes into a coproduct $\displaystyle\coprod_{\lambda \in \PP^1} \U_\lambda$,
where $\U_\lambda$ denotes the uniserial category of torsion sheaves supported on the point $\lambda$ (\cite[Proposition 2.4 and 2.5]{GL}).
In particular, Exceptional torsion sheaves can be classified as follows.
\begin{prop}[{cf. Meltzer~~\cite[Section 3.2]{M2}}]\label{prop : exceptional object}
Let the notations be as above.
For each $i=1,\dots,r$, a torsion sheaf $E \in \U_{\lambda_i}$ is exceptional if and only if the length of $E$ is smaller than $a_i$.
\qed
\end{prop}
Let $\D^b (\PP_{A,\Lambda}^1)$ be the bounded derived category of the abelian category $\coh(\PP^{1}_{A,\Lambda})$, which has several nice properties as follows.
\begin{prop}[{\cite[Section 4.1]{GL}}]\label{prop : fec can}
The ordered set of $\mu_A$ exceptional sheaves on $\PP_{A,\Lambda}^1$
\[
( \O, \O(\vec{x}_1), \cdots, \O((a_1 - 1)\vec{x}_1), \O(\vec{x}_2), \cdots, \O((a_2 - 1)\vec{x}_2),  \dots, \O(\vec{x}_r), \cdots, \O((a_r - 1)\vec{x}_r), \O(\vec{c}_A) )
\]
is a full exceptional collection in $\D^b(\PP_{A,\Lambda}^1)$.
\qed
\end{prop}
\begin{prop}[{Meltzer~\cite{M1}}]
The group ${\rm Br}_{\mu} \ltimes \ZZ^{\mu}$ act transitively on the set ${\rm FEC}(\D^b(\PP^1_{A,\Lambda}))$.
\qed
\end{prop}
\begin{prop}[{\cite[Section 2.2]{GL}}]\label{prop : Serre functor}
Define an element $\vec{\omega}_A \in L_A$ by
\begin{equation}
\vec{\omega}_A := \vec{c}_A - \sum_{i =  1}^r \vec{x}_i.
\end{equation}
Then, the autoequivalence $\S(-):= (-) \otimes \O (\vec{\omega}) [1]$ is the Serre functor of $\D^b (\PP_{A,\Lambda}^1)$.
Namely, there is an isomorphism
\[
\Hom_{\D^b(\PP_{A,\Lambda}^1)}(X, Y) \cong \Hom_{\D^b(\PP_{A,\Lambda}^1)}(Y, \S(X))^*,
\]
which is functorial with respect to $X,Y \in \D^b(\PP_{A,\Lambda}^1)$.
\qed
\end{prop}
Note that the autoequivalence $\tau:=\S(-)[-1]=(-)\otimes\O(\omega_A)$ is the Auslander--Reiten translation,
which respects the abelian category $\coh(\PP_{A,\Lambda}^1)$.
\begin{prop}[{Happel--Ringel~\cite{HaRi}}]\label{prop : AR quiver}
Suppose that $\chi_A=0$.
The components of the Auslander--Reiten quiver are of the form $\ZZ A_{\infty / r}$ for some positive integer $r\in\{1,a_1,\dots, a_r\}$.
In particular, the followings hold:
\begin{enumerate}
\item For an ordinary point $\lambda \in \PP^{1} \setminus \{ \lambda_1, \dots, \lambda_r \}$,
there is a component of the Auslander--Reiten quiver of type $\ZZ A_{\infty / 1}$ containing $S_\lambda$ as a mouth object.
\item For an orbifold point $\lambda_i$, $i=1,\dots, r$, there is a component of the Auslander--Reiten quiver of type $\ZZ A_{\infty / a_i}$
containing $S_{i;k}$ for all $k=0,\dots, a_i-1$ as mouth objects.
\qed
\end{enumerate}
Here, an indecomposable object $X$ is called a mouth object if the middle term $Y$ of the Auslander--Reiten  triangle
\[
\tau(X)\longrightarrow Y\longrightarrow X\longrightarrow \tau(X)[1]
\]
is indecomposable.
\end{prop}

Let ${\rm Pic}(\PP_{A,\Lambda}^1)$ be the {\em Picard group} of $\PP_{A,\Lambda}^1$, the subgroup of $\Aut(\coh(\PP_{A,\Lambda}^1))$
consisting of the grading shifts $(-)\otimes \O(\vec{l})$, $\vec{l}\in L_A$, which is isomorphic to $L_A$.
Denote by ${\rm Pic}_0(\PP_{A,\Lambda}^1)$ the kernel of the map $\delta: {\rm Pic}(\PP_{A,\Lambda}^1)\longrightarrow \ZZ$
defined by $(-)\otimes \O(\vec{x}_i)\mapsto \ell_A/a_i$, $i=1,\dots, r$.

Denote by $K_0(\coh(\PP_{A,\Lambda}^1))$ the \emph{Grothendieck group} of $\coh(\PP_{A,\Lambda}^1)$. It follows from \cite[Proposition 4.1]{GL} that the classes $[\O(\vec{x})] \; (0\leq \vec{x}\leq \vec{c})$ form a basis of $K_0(\coh(\PP_{A,\Lambda}^1))$. Recall from \cite{Le} that there are two important linear forms ${\rm rk}$ and ${\rm deg}$ on $K_0(\coh(\PP_{A,\Lambda}^1))$, called the \em{rank} and the \em{degree}. The rank ${\rm rk}: K_0(\coh(\PP_{A,\Lambda}^1))\rightarrow \ZZ$ is characterized by the fact that
${\rm rk}([\O(\vec{x})])=1$ for each $\vec{x}\in L_A$. The rank is zero on torsion sheaves and strictly positive on non-zero vector bundles.
Denote by the same letter $\delta$ the above homomorphism $\delta:L_A\rightarrow \ZZ$ in terms of $L_A$, which sends $\vec{x}_i$ to $\frac{\ell_A}{a_i}$ for $1\leq i\leq r$.
Then ${\rm deg}: K_0(\coh(\PP_{A,\Lambda}^1))\rightarrow \ZZ$ is characterized by the property that ${\rm deg}([\O(\vec{x})])=\delta(\vec{x}).$ It follows that $\deg([S_{\lambda}])=\ell_A$ for any $\lambda \in \PP^{1} \setminus \{ \lambda_1, \dots, \lambda_r\}$ and $\deg([S])=\frac{\ell_A}{a_i}$ if $S$ is an exceptional simple sheaf.

\begin{prop}[{Lenzing--Meltzer~\cite[Theorem~5.1, Theorem~6.3 and Proposition~7.1]{LM}}]\label{the properties of group}
Suppose $\chi_A=0$.
\begin{enumerate}
\item
There exists an isomorphism of groups
\begin{equation}
\Aut(\D^b(\PP^1_{A,\Lambda}))\cong \left({\rm Pic}_0(\PP^1_{A,\Lambda})\rtimes \Aut(\PP^1_{A,\Lambda})\right)\rtimes {\rm Br}_3.
\end{equation}
\item
The kernel of the natural group homomorphism $\pi: \Aut(\D^b(\PP^1_{A,\Lambda}))\longrightarrow \Aut(K_0(\D^b(\PP^1_{A,\Lambda})))$
is the subgroup generated by $[2]$ which will be denoted by $\ZZ[2]$.
\item
The group homomorphism $\pi$ induces the isomorphisms
\begin{equation}
{\rm Br}_3/\ZZ[2]\cong {\rm SL}(2;\ZZ),\quad {\rm Br}_3/\ZZ[1]\cong {\rm SL}(2;\ZZ)/\{\pm {\rm id}\}.
\end{equation}
where $\ZZ[1]$ denotes the subgroup generated by $[1]$.
\end{enumerate}
\qed
\end{prop}

\subsection{Congruence subroups of ${\rm SL}(2;\ZZ)$}

\begin{defn}
Let $\Gamma:={\rm SL}(2;\ZZ)$ be the special linear group of $2\times 2$ matrices with integer coefficients and determinant $1$.
Denote by $\overline{\Gamma}$ the quotient of $\Gamma$ by its center, which is $\Gamma/\{\pm {\rm id}\}$.

For each positive integer $N$, the congruence subgroup $\Gamma(N)$ is a subgroup of $\Gamma$ given by
\begin{equation}
\Gamma(N):=\left\{\left.\begin{pmatrix} a & b\\ c & d\end{pmatrix}\, \right|\, a,d\equiv 1\ ({\rm mod}\ N),\ b,c\equiv 0 \ ({\rm mod}\ N)\right\},
\end{equation}
whose image in $\overline{\Gamma}$ is denoted by $\overline{\Gamma(N)}$.
\end{defn}

It is well-known that $\overline{\Gamma}$ is the group of linear fractional transformations $\tau\mapsto (a\tau +b)/(c\tau +d)$ of the upper half plane $\HH$
and the action of $\overline{\Gamma}$ on $\HH$ is naturally extended to the one on $\overline{\HH}:=\HH\cup\{\sqrt{-1}\infty\}\cup\QQ$.
We recall the following famous facts which will be used later.
\begin{prop}[\cite{K}]
Let $N=2,3,4,6$. The index $[\overline{\Gamma}:\overline{\Gamma(N)}]$ and the number of cusps,
the equivalence class of the quotient set $(\overline{\HH}\setminus\HH)/\overline{\Gamma(N)}$, are given as follows:

\begin{table}[h]
\begin{tabular}{|c||c|c|c|c|}
\hline
$N$ & $2$ & $3$ & $4$ & $6$ \\ \hline \hline
$[\overline{\Gamma}:\overline{\Gamma(N)}]$ & $6$ & $12$ & $24$ & $72$ \\ \hline
$ (\overline{\HH}\setminus\HH)/\overline{\Gamma(N)}$ & $3$ & $4$ & $6$ & $12$\\ \hline
\end{tabular}
\vspace{5pt}
\caption{Index and the number of cusps}
\end{table}\qed
\end{prop}

\section{Statement}

Let the notations be as Section~\ref{sec : preliminaries}.
Let $\D:=\D^b(\PP_{A,\Lambda}^1)\times \D^b(\CC\vec{\Delta}_1)\times \dots \times \D^b(\CC\vec{\Delta}_m)$ for some $A$ and Dynkin quivers $\vec{\Delta}_1, \dots, \vec{\Delta}_m$.
It admits a full exceptional collection $(E_1,\dots, E_\mu)$ and the Serre functor $\S$.
The set $\FEC(\D) / \ZZ^\mu$ is of our interest but in general it is larger than what we are really want, for example, it is an infinite set.
It is natural idea to consider an action of an infinite group on it to obtain a finite number, which leads us to the subgroup of $\Aut(\D)$
generated by spherical twists.
\begin{prop}[{Seidel--Thomas~\cite[Proposition 2.10]{ST}}]
Let $S \in \D$ a spherical object, namely, $S$ is such object that
\[
\S(S) \cong S[1]\quad\text{and}\quad
\Hom_{\D} (S, S[p]) \cong
\begin{cases}
\CC, & p = 0, 1, \\
0, & p \ne 0, 1.
\end{cases}
\]
There exists an autoequivalence ${\rm Tw}_S\in \Aut(\D)$ defined by the exact triangle
\begin{equation}\label{eq : spherical twist}
\RR\Hom_{\D} (S, X) \otimes S \longrightarrow X \longrightarrow {\rm Tw}_S(X)\longrightarrow (\RR\Hom_{\D} (S, X) \otimes S )[1],\quad X\in \D.
\end{equation}
The inverse functor ${\rm Tw}_S^{-1}\in \Aut(\D)$ is given by
\begin{equation}\label{eq : spherical twist}
{\rm Tw}_S^{-1}(X)\longrightarrow X\longrightarrow S\otimes\RR\Hom_{\D} (X, S)^{*}\longrightarrow ({\rm Tw}_S^{-1}(X))[1].
\end{equation}
\qed
\end{prop}

We define a subgroup ${\rm ST}(\D)$ of the group of autoequivalences $\Aut(\D)$,
which stands for the initial of Seidel--Thomas or spherical twists.
\begin{defn}
Define a group ${\rm ST}(\D)$ as the set generated by spherical twists:
\begin{equation}
{\rm ST}(\D) \coloneqq \big\langle {\rm Tw}_S \in \Aut(\D) \mid S \in {\rm Sph}(\D) \big\rangle,
\end{equation}
where ${\rm Sph}(\D)$ is the set of spherical objects in $\D$.
\end{defn}

The reason we use ${\rm Sph}(\D)$ is since we are mainly interested in triangulated categories
``associated to/mirror dual to some geometry of algebraic curves/Riemannian surfaces".
Indeed, the group ${\rm ST}(\D)$ has a closely related property to the mapping class group of a Riemannian surface,
which we call the {\em Alexander method}.
\begin{lem}[Alexander method]\label{lem:alex}
Let $E \in \D$ be an exceptional object.
For a spherical object $S \in {\rm Sph}(\D)$, ${\rm Tw}_S (E) \cong E$ if and only if $S\in E^\perp$.
In particular, we have
\[
{\rm ST}(E^\perp) \cong \big\langle {\rm Tw}_S \in {\rm ST}(\D) \mid S \in {\rm Sph}(\D), ~ {\rm Tw}_S (E) \cong E \big\rangle .
\]
\qed
\end{lem}
Here, due to \cite[Theorem~3.3.2 and Theorem~3.3.3]{M2}, $E^\perp$ is of the form
$\D^b(\PP_{A'}^1)\times \D^b(\CC\vec{\Delta}_1)\times \dots \times \D^b(\CC\vec{\Delta}_{m'})$ for some $A'$ and Dynkin quivers $\vec{\Delta}_1, \dots, \vec{\Delta}_{m'}$.

Thus, in view of Lemma~\ref{lem:alex}, we define $e(\D)$ by
\begin{equation}
e(\D) := |\FEC(\D) / \langle {\rm ST}(\D) , \ZZ^{\mu}\rangle |.
\end{equation}

\begin{prop}[{\cite{De, ONSFR}}]\label{prop : number of full exceptional collections for Dynkin}
Let $\vec{\Delta}$ be a Dynkin quiver. Let $\D^b(\CC\vec{\Delta})$ be the bounded derived category of finite dimensional modules over the path algebra $\CC\vec{\Delta}$
and $\mu:={\rm rank}_\ZZ K_0(\D^b(\CC\vec{\Delta}))$.
We have ${\rm ST}(\D^b(\CC\vec{\Delta})) = \{ 1 \}$ and
\begin{equation}
|\FEC(\D^b(\CC\vec{\Delta})) / \ZZ^\mu|= \frac{\mu!}{d_1\cdots d_\mu}h^\mu.
\end{equation}
where $2\le d_1\le\dots\le d_\mu=h$ are degrees of $\mu$ algebraic independent invariants of the Weyl group for the Dynkin diagram $\Delta$.
\begin{table}[h]
\centering
\begin{tabular}{|c||c|c|c|c|c|} \hline
$\Delta$ & $A_\mu$ & $D_\mu$ & $E_6$ & $E_7$ & $E_8$ \\ \hline
$e(\D^b(\CC\vd))$ & $(\mu + 1)^{\mu - 1}$ & $2(\mu - 1)^\mu$ & $2^9 \cdot 3^4$ & $2 \cdot 3^{12}$ & $2 \cdot 3^5 \cdot 5^7$ \\ \hline
\end{tabular}
\vspace{5pt}
\caption{Number of full exceptional collections for Dynkin quivers}
\end{table}
\qed
\end{prop}

\begin{prop}[{Otani--Shiraishi--Takahashi~\cite[Proposition~4.8 and Theorem~4.10]{OST}}]
If $\chi_A>0$, then ${\rm ST}(\D^b(\PP^1_{A,\Lambda}))=\langle (-)\otimes\O(\vec{c})\rangle$.
In particular, we have
\begin{equation}
e(\D^b(\PP_{A,\Lambda}^1)) =\frac{\mu_A!}{a_1!\cdots a_r!\chi_{A}}a_1^{a_1}\cdots a_r^{a_r}.
\end{equation}\qed
\end{prop}

Consider the quotient of the subgroup ${\rm ST}(\D^b(\PP^1_{A,\Lambda}))$ of $\Aut(\D^b(\PP^1_{A,\Lambda}))$ by ${\rm ST}(\D^b(\PP^1_{A,\Lambda}))\cap \ZZ[1]$
and denote it by $\overline{{\rm ST}(\D^b(\PP^1_{A,\Lambda}))}$.
Since the abelian category ${\rm coh}\,\PP^1_{A,\Lambda}$ is hereditary, we have
\[
e(\D^b(\PP_{A,\Lambda}^1)) =| \{(E_1,\dots, E_{\mu_A})\in \FEC(\D^b(\PP_{A,\Lambda}^1))\,|\, E_1,\dots, E_{\mu_A}\in \coh(\PP_{A,\Lambda}^1)\} /\overline{{\rm ST}(\D^b(\PP^1_{A,\Lambda}))}|.
\]

\begin{thm}\label{thm:main1}
Suppose $\chi_A=0$. The group ${\rm ST}(\D^b(\PP^1_{A,\Lambda}))$ is a subgroup of ${\rm Br}_3$ satisfying
\begin{equation}
{\rm ST}(\D^b(\PP^1_{A,\Lambda}))/\ZZ[2]\cong \Gamma(\ell_A),\quad \overline{{\rm ST}(\D^b(\PP^1_{A,\Lambda}))}\cong \overline{\Gamma(\ell_A)}.
\end{equation}
\end{thm}
\begin{thm}\label{thm:main2}
Suppose $\chi_A=0$.
Then we have
\begin{equation}
e(\D^b(\PP_{A,\Lambda}^1)) =\frac{\mu_A!}{a_1!\cdots a_r!}a_1^{a_1}\cdots a_r^{a_r}\cdot
\frac{\displaystyle\sum_{i=1}^r a_i^2 (a_i-1)}{2\mu_A}\cdot \frac{[\overline{\Gamma}:\overline{\Gamma(\ell_A)}]}{\ell_A}.
\end{equation}
In particular,
\[
|{\rm FEC}(\D^b(\PP_{A,\Lambda}^1))/\overline{\Gamma}|=
\frac{e(\D^b(\PP_{A,\Lambda}^1))}{[\overline{\Gamma}:\overline{\Gamma(\ell_A)}]}
=
\frac{\mu_A!}{a_1!\cdots a_r!}a_1^{a_1}\cdots a_r^{a_r}\cdot
\frac{\displaystyle\sum_{i=1}^r a_i^2 (a_i-1)}{2\ell_A\mu_A}.
\]
\end{thm}

The following gives a consistency with the expectation from the mirror symmetry of $\PP^1_{A,\Lambda}$ that
${\rm Stab}(\D^b(\PP^1_{A,\Lambda}))$ is connected and simply connected and its quotient by $\overline{\Gamma(2)}$ can be identified with the base space $\CC^{\mu_A-1}\times (\CC\setminus\{0,1\})$ of the universal unfolding of its mirror dual.
Note here that $\overline{\Gamma(2)}$ is isomorphic to the free group $F_2$ of rank two, which is the fundamental group of $\CC\setminus\{0,1\}$.
\begin{cor}\label{cor:3.8}
Suppose that $A=(3,3,3), (2,4,4)$ or $(2,3,6)$. The number
\[
|{\rm FEC}(\D^b(\PP_{A,\Lambda}^1))/\overline{\Gamma(2)}|=
\frac{e(\D^b(\PP_{A,\Lambda}^1))}{[\overline{\Gamma}:\overline{\Gamma(\ell_A)}]}\cdot [\overline{\Gamma}:\overline{\Gamma(2)}]
=\frac{\mu_A!}{a_1!\cdots a_r!}a_1^{a_1}\cdots a_r^{a_r}\cdot
\frac{\displaystyle\sum_{i=1}^r a_i^2 (a_i-1)}{2\ell_A\mu_A}\cdot 6,
\]
coincides with the degree ${\rm deg} LL^{alg}$ of the Lyashko--Looijenga map by Hertling--Roucairol~\cite[Theorem~6.3]{HR},
for the universal unfolding of the simple elliptic singularities in the Legendre normal forms with $A$ as Gabrielov numbers (cf. Ebeling--Takahashi~\cite{ET}):
\begin{table}[h]
\begin{tabular}{|c||c|c|c|}
\hline
Singularity Type  & $A$  & Legendre Normal Form $f$ \\ \hline \hline
$\widetilde{E}_6$ & $(3,3,3)$ & $z_2(z_2-z_1)(z_2-\lambda z_1)-z_1z_3^2$\\ \hline
$\widetilde{E}_7$ & $(2,4,4)$ & $z_1z_2(z_2-z_1)(z_2-\lambda z_1)+z_3^2$\\ \hline
$\widetilde{E}_8$ & $(2,3,6)$ & $z_2(z_2-z_1^2)(z_2-\lambda z_1^2)+z_3^2$\\ \hline
\end{tabular}
\vspace{5pt}
\caption{Gabrielov numbers and the Legendre normal forms}
\end{table}
\qed
\end{cor}
\begin{rem}
For $A=(2,2,2,2)$, we have $e(\PP^1_{A,\Lambda})=6!\cdot 4\cdot 2^4$, which is a Hurwitz number (of genus $0$ with $4$ boundaries)
$6!\cdot 4$ times $2^4$ and can also be considered as the degree of a Lyashko--Looijenga map of the Frobenius manifold studied by Basalaev~\cite{Ba}.
On the relation of number of full exceptional collections with Hurwitz numbers, we will make it clear elsewhere in the future.
\end{rem}

\section{Proof of Theorem~\ref{thm:main1}}
Following \cite[Section~6]{LM} we take the subgroup ${\rm Br}_3$ of $\Aut(\D^b(\PP_{A,\Lambda}^1))$ so that we have
\[
{\rm Br}_3=\langle b_1,b_2\,|\, b_1b_2b_1=b_2b_1b_2\rangle,\quad b_1(-)=(-)\otimes\O(\vec{x_r}),
\]
where $A=(a_1,\dots ,a_r)$ is chosen as in Table~\ref{tab: classification A},
and the group homomorphism $\pi: {\rm Br}_3\longrightarrow \Gamma$, induced by the group homomorphism
$\pi:\Aut(\D^b(\PP_{A,\Lambda}^1))\longrightarrow \Aut(K_0(\D^b(\PP_{A,\Lambda}^1)))$, is given as
\[
b_1\mapsto \begin{pmatrix} 1& 1\\ 0 & 1\end{pmatrix},\quad b_2\mapsto \begin{pmatrix} 1& 0\\ -1 & 1\end{pmatrix}.
\]
\begin{lem}\label{the property of spherical twist}
The spherical twist ${\rm Tw}_{S_\lambda}$ for a simple sheaf $S_\lambda$, $\lambda \in \PP^{1} \setminus \{ \lambda_1, \dots, \lambda_r \}$
is isomorphic to the autoequivalence $(-)\otimes\O(\vec{c}_A)\in\Aut(\D^b(\PP_{A,\Lambda}^1))$.
In particular, we have ${\rm Tw}_{S_\lambda}=b_1^{\ell_A}\in {\rm Br}_3$.
\end{lem}
\begin{pf}
Since $S_\lambda\otimes \O(\vec{l})\cong S_\lambda$ for $\vec{l}\in L_A$, the short exact sequence \eqref{eq:ordinary simple} yields the following exact triangle
\[
S_\lambda[-1] \longrightarrow \O(\vec{l})\xrightarrow{x_2^{a_2}-\lambda x_1^{a_1}}\O(\vec{l}+\vec{c}_A)\longrightarrow S_\lambda,\quad \vec{l}\in L_A,
\]
which implies ${\rm Tw}_{S_\lambda}(-)=(-)\otimes\O(\vec{c}_A)$ due to Proposition~\ref{prop : fec can}.
Since $\vec{c}_A=a_r\vec{x}_r=\ell_A\vec{x}_r$, the statement follows.
\qed
\end{pf}
Note that  elements of $\Aut(\PP^1_{A,\Lambda})$ send ordinary points to ordinary points. Thus we have the following
\begin{lem}
For $S\in{\rm Sph}(\D^b(\PP_{A,\Lambda}^1))$,
there exists $b\in {\rm Br}_3$ such that $b(S)\cong S_\lambda$ for some ordinary point $\lambda \in \PP^{1} \setminus \{ \lambda_1, \dots, \lambda_r \}$.
In particular, ${\rm Tw}_S=b^{-1}{\rm Tw}_{b(S)}b=b^{-1}b_1^{\ell_A} b\in {\rm Br}_3$.
\end{lem}
\begin{pf}
According to \cite{Le}, any indecomposable object belongs to $\coh(\PP_{A,\Lambda}^1)[n]$ for some $n$.
Therefore, combining with Proposition \ref{the properties of group} (iii), we may assume that $S\in \coh(\PP_{A,\Lambda}^1)$ or $S\in\coh(\PP_{A,\Lambda}^1)[1]$.
Note that such object $S$ are characterized by the pair of integers $({\rm deg}(S), {\rm rk}(S))\in\ZZ^2$ and a point in $\PP^1$.

Since $S$ is a spherical object, it follows that $S$ is a mouth object in the homogeneous tube. Combining ${\rm gcd}({\rm deg}(S_{\mu}), {\rm rk}(S_{\mu}))=(\ell_A, 0)=\ell_A$ for any $\mu \in \PP^{1} \setminus \{ \lambda_1, \dots, \lambda_r\}$ with the exact sequence (\ref{eq:ordinary simple}), we obtain that ${\rm gcd}({\rm deg}(S), {\rm rk}(S))=\ell_A$. On the other hand, there exists $A\in {\rm SL}(2;\ZZ)$ such that $({\rm deg}(S), {\rm rk}(S))A=(\ell_A, 0).$ By Kac's Theorem \cite{CB} we obtain that an indecomposable sheaf $S'$ satisfies $({\rm deg}(S'), {\rm rk}(S'))=(\ell_A, 0)$ if and only if $S'$ is isomorphic to some ordinary simple sheaf $S_\lambda,$ which yields the first statement.

Combining the Lemma \ref{the property of spherical twist} and the action of braid group is transitive due to \cite{M1}, we have $T_S=b^{-1}T_{b(S)}b=b^{-1}b_1^{\ell_A} b\in {\rm Br}_3$.
\qed
\end{pf}
Hence, it follows that ${\rm ST}(\D^b(\PP_{A,\Lambda}^1))$ is a subgroup of ${\rm Br}_3$.
Moreover, we have
\[
{\rm ST}(\D^b(\PP^1_{A,\Lambda}))/\ZZ[2]\stackrel{\pi}{\cong} \left\langle\left.M^{-1}\cdot \begin{pmatrix} 1& \ell_A\\ 0 & 1\end{pmatrix}\cdot M\,\right|\, M\in \Gamma\right\rangle,
\]
which is obviously a normal subgroup of $\Gamma$.
It is a subgroup of $\Gamma(\ell_A)$ since $\Gamma(\ell_A)$ is a normal subgroup of $\Gamma$.
On the other hand, by a direct calculation on a case-by-case basis  it turns out that each generator $g$ of $\Gamma(\ell_A)$ belongs to the RHS of the above isomorphism.
Therefore, we have ${\rm ST}(\D^b(\PP^1_{A,\Lambda}))/\ZZ[2]\cong \Gamma(\ell_A)$.
The rest of the statement is clear.
\section{Proof of Theorem~\ref{thm:main2}}

For each exceptional simple sheaf $S_{i;k}$ there is a unique exceptional sheaf $E_{i,j;k}$ which is uniserial of length $j=1,\dots, a_i-1$ such that
\[
E_{i,a_i-j;k}\longrightarrow E_{i,j+1;k}\longrightarrow \dots\longrightarrow E_{i,a_i-1;k}=S_{i;k}
\]
is a sub-quiver of the Auslander--Reiten quiver.

Suppose that the last term $E_{\mu_A}$ of a full exceptional collection $(E_1,\dots, E_{\mu_A})\in {\rm FEC}(\D^b(\PP^1_{A,\Lambda}))$ is $E_{i,j;k}$.
Then the semi-orthogonal complement $E_{\mu_A}^\perp=(E_{i,j;k})^\perp$ is equivalent to
\[
\D^b({\rm coh}\,\PP^1_{A_{i,j},\Lambda})\times \D^b({\rm mod}(\CC\vec{A}_{j-1})),
\]
where $A_{i,j}:=(a_1,\dots a_{i-1},a_i-j,a_{i+1}\dots ,a_r)$ and we regard $\D^b({\rm mod}(\CC\vec{A}_{0}))$ as the zero category.

For an exceptional sheaf $E$ in ${\rm coh}\,\PP^1_{A,\Lambda}$, there exists an element $\bar{\gamma}\in Br_3/\ZZ[1]=\overline{\Gamma}$
such that $\bar{\gamma}(E)=E_{i,j;k}$ for some $i,j,k$.
The above description of the semi-orthogonal complement together with Lemma~\ref{lem:alex} gives
a correspondence between exceptional objects and the cusps of the group $\overline{\Gamma(\ell_A)}$.
Namely, since the subgroup $\langle (-)\otimes\O(\vec{c})\rangle\subseteq {\rm ST}(\D^b(\PP^1_{A,\Lambda}))$ acts trivially on $E_{i,j;k}$ but
defines a nontrivial subgroup of $\overline{\Gamma}$ due to Lemma~\ref{lem:alex}, we see that an infinite cyclic subgroup of $\overline{\Gamma}$ acts trivially on $E$.
Note that for $\overline{\Gamma(\ell_A)}$ the number of cusps are $[\overline{\Gamma}:\overline{\Gamma(\ell_A)}]/\ell_A$.
Therefore, we obtain the following recursion relation:
\begin{equation}\label{eq:recursion}
e(\D^b(\PP^1_{A,\Lambda}))=\sum_{i=1}^r\sum_{j=1}^{a_i-1}\binom{\mu_A-1}{j-1}\cdot e(\D^b(\PP^1_{A_{i,j},\Lambda}))\cdot e(\D^b(\CC\vec{A}_{j-1}))\cdot a_i\cdot \frac{[\overline{\Gamma}:\overline{\Gamma(\ell_A)}]}{\ell_A}.
\end{equation}

For simplicity, we consider the part of the above sum with $i=1$. Recall that
\[
e(\D^b(\PP^1_{A_{1,j},\Lambda}))=\frac{(\mu_A-j)!}{(a_1-j)!a_2!\cdots a_r!\chi_{A_{1,j}}}(a_1-j)^{a_1-j}a_2^{a_2}\cdots a_r^{a_r},\quad e(\D^b(\CC\vec{A}_{j-1}))=j^{j-2}.
\]
In particular, since $\chi_A=0$ it turns out that
\[
\chi_{A_{1,j}}=2+\frac{1}{a_1-j}-1+\sum_{i=2}^r\left(\frac{1}{a_i}-1\right)=
\frac{1}{a_1-j}-\frac{1}{a_1}+\chi_A=\frac{j}{a_1(a_1-j)}.
\]

\begin{lem}
We have the following identities:
\begin{subequations}
\begin{equation}\label{ab:1}
nx^{n-1}=\sum_{j=1}^n\binom{n}{j}(x-j)^{n-j}j^{j-1},
\end{equation}
\begin{equation}\label{ab:2}
\frac{n(n-1)}{2}x^{n-2}=\sum_{j=1}^n\binom{n}{j}(x-j)^{n-j}(j-1)j^{j-2}.
\end{equation}
\end{subequations}
\end{lem}
\begin{pf}
It follows from Abel's identity
\[
(x+y)^n=\sum_{j=0}^n\binom{n}{j}(x-j)^{n-j}y(y+j)^{j-1}
\]
that
\[
\frac{(x+y)^n-x^n}{y}=\sum_{j=1}^n\binom{n}{j}(x-j)^{n-j}(y+j)^{j-1}.
\]
Setting $y=0$, we obtain the first identity.
Taking the partial derivative of both sides with respect to $y$ and then setting $y=0$, we obtain the second one.
\qed
\end{pf}
\begin{cor}
We have
\begin{equation}
\frac{1}{2}(n-1)n^{n}=\sum_{j=1}^n\binom{n}{j}(n-j)^{n-j+1}j^{j-2}.
\end{equation}
\end{cor}
\begin{pf}
It follows from the identity \eqref{ab:1} and \eqref{ab:2} that
\begin{eqnarray*}
\frac{n(n-1)}{2}x^{n-1}&=&x\sum_{j=1}^n\binom{n}{j}(x-j)^{n-j}(j-1)j^{j-2}\\
&=&\sum_{j=1}^n\binom{n}{j}(x-j)^{n-j}(xj-j-x+j)j^{j-2}\\
&=&nx^{n}-nx^{n-1}-\sum_{j=1}^n\binom{n}{j}(x-j)^{n-j+1}j^{j-2}.
\end{eqnarray*}
Setting $x=n$, we obtain the statement.
\qed
\end{pf}

\begin{lem}
We have
\begin{eqnarray*}
& &\sum_{j=1}^{a_1-1}\binom{\mu_A-1}{j-1}\frac{(\mu_A-j)!}{(a_1-j)!a_2!\cdots a_r!\frac{j}{a_1(a_1-j)}}(a_1-j)^{a_1-j}a_2^{a_2}\cdots a_r^{a_r} j^{j-2}\\
&=&\frac{\mu_A!}{a_1!a_2!\cdots a_r!}a_1^{a_1}a_2^{a_2}\cdots a_r^{a_r}\cdot \frac{a_1(a_1-1)}{2\mu_A}.
\end{eqnarray*}
\end{lem}
\begin{pf}
Applying the previous Corollary, it follows that
\begin{eqnarray*}
& &\sum_{j=1}^{a_1-1}\binom{\mu_A-1}{j-1}\frac{(\mu_A-j)!}{(a_1-j)!a_2!\cdots a_r!\frac{j}{a_1(a_1-j)}}(a_1-j)^{a_1-j}a_2^{a_2}\cdots a_r^{a_r} j^{j-2}\\
&=&\frac{\mu_A!}{a_1!a_2!\cdots a_r!}a_1^{a_1}a_2^{a_2}\cdots a_r^{a_r}\cdot \frac{1}{\mu_A}\cdot \sum_{j=1}^{a_1-1}\binom{a_1}{j}\frac{(a_1-j)^{a_1-j+1}j^{j-2}}{a_1^{a_1-1}}\\
&=&\frac{\mu_A!}{a_1!a_2!\cdots a_r!}a_1^{a_1}a_2^{a_2}\cdots a_r^{a_r}\cdot \frac{a_1(a_1-1)}{2\mu_A}.
\end{eqnarray*}
\qed
\end{pf}

The calculation is similar for $i=2,\dots, r$. Therefore, we can calculate the RHS of \eqref{eq:recursion} as
\begin{eqnarray*}
e(\D^b(\PP^1_{A,\Lambda}))&=&\sum_{i=1}^r\sum_{j=1}^{a_i-1}\binom{\mu_A-1}{j-1}\cdot e(\D^b(\PP^1_{A_{i,j},\Lambda}))\cdot e(\D^b(\CC\vec{A}_{j-1}))\cdot a_i\cdot \frac{[\overline{\Gamma}:\overline{\Gamma(\ell_A)}]}{\ell_A}\\
&=&\frac{\mu_A!}{a_1!\cdots a_r!}a_1^{a_1}\cdots a_r^{a_r}\cdot \sum_{i=1}^r\left(\frac{a_i(a_i-1)}{2\mu_A}\cdot a_i\cdot \frac{[\overline{\Gamma}:\overline{\Gamma(\ell_A)}]}{\ell_A}\right),\\
&=&\frac{\mu_A!}{a_1!\cdots a_r!}a_1^{a_1}\cdots a_r^{a_r}\cdot \frac{\displaystyle\sum_{i=1}^ra_i^2(a_i-1)}{2\mu_A}\frac{[\overline{\Gamma}:\overline{\Gamma(\ell_A)}]}{\ell_A}.
\end{eqnarray*}
We have finished the proof of Theorem~\ref{thm:main2}.


\end{document}